\documentclass[12pt]{article}
\usepackage{amsmath}
\usepackage{amssymb}
\usepackage{latexsym}
\usepackage{amssymb}
\usepackage{epsf}

\newtheorem{Theorem}{Theorem}
\newtheorem{Lemma}{Lemma}
\newtheorem{Proposition}{Proposition}

\numberwithin{equation}{section}

\begin{document}
\date{}
\author{M.I.Belishev \thanks {Saint-Petersburg Department of
                 the Steklov Mathematical Institute (POMI), 27 Fontanka,
                 St. Petersburg 191023, Russia; belishev@pdmi.ras.ru.
                 Supported by RFBR grants 11-01-00407A and NSh-4210.2010.1.}}

\title{Dynamical System with Boundary Control Associated with Symmetric
Semi-Bounded Operator}

\maketitle
\begin{abstract}
Let $L_0$ be a closed densely defined symmetric semi-bounded
operator with nonzero defect indexes in a separable Hilbert space
$\cal H$. It determines a {\it Green system} $\{{\cal H}, {\cal
B}; L_0, \Gamma_1, \Gamma_2\}$, where ${\cal B}$ is a Hilbert
space, and $\Gamma_i: {\cal H} \to \cal B$ are the operators
related through the Green formula $$(L_0^*u, v)_{\cal
H}-(u,L_0^*v)_{\cal H}=(\Gamma_1 u, \Gamma_2 v)_{\cal B} -
(\Gamma_2 u, \Gamma_1 v)_{\cal B}.$$ The {\it boundary operators}
$\Gamma_i$ are chosen canonically in the framework of the Vishik
theory.

With the Green system one associates a {\it dynamical system with
boundary control} (DSBC)
\begin{align*}
& u_{tt}+L_0^*u = 0  && {\rm in}\,\,\,{\cal H}, \,\,\,t>0\\
& u|_{t=0}=u_t|_{t=0}=0 && {\rm in}\,\,\,{\cal H}\\
& \Gamma_1 u = f && {\rm in}\,\,\,{\cal B},\,\,\,t \geqslant 0.
\end{align*}
We show that this system is {\it controllable} if and only if the
operator $L_0$ is completely non-self-adjoint.

A version of the notion of a {\it wave spectrum} of $L_0$ is
introduced. It is a topological space determined by $L_0$ and
constructed from reachable sets of the DSBC.
\end{abstract}
\setcounter{section}{-1}

\section{Introduction}
\subsection{About the paper}
We develop ideas and results of the papers \cite{DSBC} and
\cite{BArXiv}. The future prospect and goal is a {\it functional
model} of a symmetric semi-bounded operator outlined in
\cite{BArXiv}. Our paper is a step towards this model, which
prepares two of the model basic elements: {\it DSBC} and {\it wave
spectrum}.

Motivation comes from inverse problems. Namely, the inspiring role
is played by the problem of reconstruction of a Riemannian
manifold via its boun\-dary inverse data
\cite{BIP97},\cite{BIP07},\cite{BSobolev}. In accordance with the
program, which we promote, to solve the latter problem (and a
class of closely related problems) is to construct a certain
functional model of a relevant symmetric operator. The operator
codes information about the manifold and is determined by inverse
data. Its wave spectrum turns out to be isometric to the manifold.
By the latter, to decode the information one can find the wave
spectrum \cite{BArXiv}.

The main subjects of the paper are the following.

\subsection{Operator $L_0$}
Let $\cal H$ be a (separable) Hilbert space, $L_0$ a closed
operator in $\cal H$ such that ${\rm clos\,\,Dom\,}L_0=\cal H$,
$L_0\subset L_0^*$, $(L_0y,y)\geqslant \varkappa
\|y\|^2\,\,\,(\varkappa=\rm const)$. In what follows, without loss
of generality, we deal with $\varkappa>0$.

Also, we assume that $L_0$ has nonzero defect indexes:
$n_+=n_-=\newline{\rm dim\,\,Ker\,}L_0^* \geqslant 1$. Note that
such an operator is necessarily unbounded.

No more assumptions on $L_0$ are imposed. It is a class of
operators, for which we plan to construct the above-mentioned
functional model.

\subsection{Green system}
Let $L$ be the {\it Friedrichs extension} of $L_0$, so that $L_0
\subset L \subset L_0^*,\,\,\,L^*=L$, and $(L_0y,y)\geqslant
\varkappa \|y\|^2$ holds (see, e.g., \cite{BirSol}). The inverse
$L^{-1}$ is bounded and defined on $\cal H$.

The well-known decomposition by Vishik \cite{Vishik} is $${\rm
Dom\,}L_0^*={\rm Dom\,}L_0 \overset{.}+L^{-1}{\rm
Ker\,}L_0^*\overset{.}+{\rm Ker\,}L_0^*$$ (direct sums). By this,
for a $y \in {\rm Dom\,}L_0^*$ one has $$y=y_0+L^{-1}g+h$$ with
$y_0 \in {\rm Dom\,}L_0$ and $g,h \in {\rm Ker\,}L_0^*$. Denote
$h=\Gamma_1y$ and $g=\Gamma_2y$. We derive the {\it Green formula}
$$(L_0^*u, v)_{\cal H}-(u,L_0^*v)_{\cal H}=(\Gamma_1 u, \Gamma_2
v)_{{\rm Ker\,}L_0^*} - (\Gamma_2 u, \Gamma_1 v)_{{\rm
Ker\,}L_0^*}\,,$$ which is in fact a partial case of more general
relation established in \cite{Vishik}. Hence, a collection
$\{{\cal H}, {\rm Ker\,}L_0^*; L_0, \Gamma_1, \Gamma_2\}$
constitutes a {\it Green system} determined by the operator $L_0$.

\subsection{DSBC}
The Green system, in turn, determines a {\it dynamical system with
boundary control}
\begin{align*}
& u_{tt}+L_0^*u = 0  && {\rm in}\,\,\,{\cal H}, \,\,\,t>0\\
& u|_{t=0}=u_t|_{t=0}=0 && {\rm in}\,\,\,{\cal H}\\
& \Gamma_1 u = f && {\rm in}\,\,\,{{\rm Ker\,}L_0^*},\,\,\,t
\geqslant 0,
\end{align*}
where $f=f(t)$ is a ${\rm Ker\,}L_0^*$-valued function of time
({\it boundary control}). By $u=u^f(t)$ we denote the
(generalized) solution, which is well defined for a class $\cal M$
of smooth enough controls $f$.

A set $${\cal U}^t:=\{u^f(t)\,|\,\,f \in \cal M\}$$ is called {\it
reachable} (at the moment $t$), whereas $${\cal U}\,:=\,\bigvee_{t
\geqslant 0}{\cal U}^t$$ (algebraic sum) is a {\it total reachable
set}. The DSBC is said to be {\it controllable} if
\begin{equation*}
{\rm clos\,}{\cal U}\,=\,\cal H\,.
\end{equation*}
We prove that this relation holds if and only if $L_0$ is {\it
completely non-self-adjoint operator}. The latter means that there
is no nonzero subspace in $\cal H$, in which $L_0$ has a
self-adjoint part.

\subsection{Wave spectrum}
This notion is introduced in a few steps.
\smallskip

First, we define a so-called {\it inflation} ${I_{L_0}}$, which is
an operation on the lattice ${{\mathfrak L}(\cal H)}$ of subspaces
in $\cal H$. Inflation extends subspaces and is determined by the
operator $L_0$ (more precisely, by its Friedrichs extension $L$).

Second, we define ${{\mathfrak L}_{L_0}}$ as the {\it minimal
sublattice} of ${{\mathfrak L}(\cal H)}$, which contains all
reachable subspaces ${\rm clos\,}{\cal U}^t,\,\,t \geqslant 0$ and
is invariant with respect to the inflation ${I_{L_0}}$.

Next, we introduce a family ${I_{L_0}}{{\mathfrak L}_{L_0}}$ of
monotone (growing) ${{\mathfrak L}_{L_0}}$-valued functions of $t
\geqslant 0$ provided with standard lattice topology. This family
is a partially ordered set. As such, it may content a set ${\rm
At\,}{I_{L_0}}{{\mathfrak L}_{L_0}}$ of minimal nonzero elements
({\it atoms}) .

At last, ${\rm At\,}{I_{L_0}}{{\mathfrak L}_{L_0}}$ is endowed
with a relevant ({\it ball}\,-) topology $\beta$, and we arrive at
a topological space $({\Omega_{L_0}}, \beta)$. It is the space,
which we call a {\it wave spectrum} of the operator $L_0$.

\subsection{Reconstruction of manifolds}
As was noted in 0.1, our program is motivated by inverse problems.
As is shown in \cite{BArXiv}, to recover a Riemannian manifold
$\Omega$ from the boundary inverse data one can
\begin{itemize}
\item determine a unitary copy $\tilde L_0$ of the {\it minimal
Laplacian} $L_0=-\Delta$ in $\Omega$ from the data \item find the
wave spectrum $\Omega_{\tilde L_0}$.
\end{itemize}
For a generic class of manifolds, the space $\Omega_{\tilde L_0}$
turns out to be isometric to $\Omega$. By this, the wave spectrum
provides a representative of the class of manifolds, which possess
the given inverse data. Thus, $\Omega_{\tilde L_0}$ solves the
reconstruction problem.

\section{DSBC}

\subsection{Operator $L_0$}
Let us specify the class of operators, which we deal with.

Let $\cal H$ be a (separable) Hilbert space, $L_0$ an operator in
$\cal H$. We assume that
\begin{enumerate}
\item $L_0$ is closed and densely defined: ${\rm clos\,}{\rm
Dom\,}L_0 = \cal H$ \item $L_0$ is positive definite: there is a
constant $\varkappa
>0$ such that $\\(L_0 y,y)\geqslant \varkappa \|y\|^2$ holds as $y \in {\rm Dom\,}L_0$
\item $L_0$ has the nonzero defect indexes $n_+=n_-=\dim {\rm
Ker\,} L_0^*$: $\\1 \leqslant \dim {\rm Ker\,}L_0^* \leqslant
\infty.$
\end{enumerate}
Note that by 3. such an operator is necessarily unbounded.
\smallskip

By $L$ we denote the Friedrichs extension of $L$, so that $L_0
\subset L \subset L_0^*,\,\,L^*= L,$ and $(L y,y)\geqslant
\varkappa \|y\|^2$ hold for all $y \in {\rm Dom\,}L$ (see, e.g.,
\cite{BirSol}). Its inverse $L^{-1}$ is a bounded operator defined
on $\cal H$.

\subsection{Green system}
We begin with basic definitions, which go back to the pioneer
paper by A.N.Kochubei \cite{Koch} (see also \cite{Ryzh}).
\smallskip

Let $\cal H$ and $\cal B$ be the Hilbert spaces, $A: {\cal H} \to
\cal H$ and $\Gamma_i: {\cal H} \to {\cal B}\,\,\,(i=1,2)$ the
operators such that ${\rm clos\,\,Dom\,}A={\cal H},\,\,{\rm Dom\,}
\Gamma_i \supset {\rm Dom\,}A,\,\,{\rm clos\,} \vee_{i=1,2}{\rm
Ran\,}\Gamma_i = \cal B $.

A collection ${\mathfrak G}\!{\mathfrak r}=\{{\cal H}, {\cal B};
A, \Gamma_1, \Gamma_2\}$ is said to be a {\it Green system}, if
its elements are related via the {\it Green formula}
\begin{equation}\label{Green Formula}
(Au,v)_{\cal H}-(u,Av)_{\cal H}=(\Gamma_1u, \Gamma_2 v)_{\cal B} -
(\Gamma_2u, \Gamma_1 v)_{\cal B}
\end{equation}
for all $u, v \in {\rm Dom\,} A$. The space $\cal H$ is called
{\it inner}, $\cal B$ is a {\it space of boundary values}, $A$ is
a {\it basic operator}, $\Gamma_{0, 1}$ are the {\it boundary
operators}.
\medskip

In a Green system with the given $\cal H, \cal B, A$, there is a
freedom of choice of the boundary operators. For instance, taking
an $S=S^*, \,\,{\rm Dom\,}S \supset {\rm Ran\,}\Gamma_1$ and
putting $\widetilde \Gamma_2:= \Gamma_2 + S\Gamma_1$, one gets a
collection $\widetilde {{\mathfrak G}\!{\mathfrak r}}= \{{\cal H},
{\cal B}; A, \Gamma_1, \widetilde \Gamma_2\}$, which is also a
Green system.

\subsection{System ${{\mathfrak G}\!{\mathfrak r}}_{L_0}$}
Here we associate with $L_0$ a Green system with the canonically
chosen boundary operators.
\smallskip

Denote $${\cal K}\,:=\,{\rm Ker\,} L_0^*$$ and recall that $\dim
{\cal K} \geqslant 1$. Let $P$ be the (orthogonal) projection in
$\cal H$ onto ${\cal K}$, $\mathbb O$ and $\mathbb I$ the zero and
unit operators. Also, introduce the operators
$$\Gamma_1\,:=\,L^{-1}L_0^*-{\mathbb I}\,, \quad \Gamma_2\,:=\,P L_0^*.$$
\begin{Lemma}\label{Lemma 1}
The collection ${{{\mathfrak G}{\mathfrak r}}_{L_0}}:=\{{\cal H},
{\cal K};\,L_0^*, \Gamma_1, \Gamma_2\}$ is a Green system.
\end{Lemma}
{\bf Proof}
\smallskip

\noindent {\bf 1.}\,\,\,Since $L_0^*
\Gamma_1=L_0^*L^{-1}L_0^*-L_0^*=L L^{-1}L_0^*-L_0^*=\mathbb O$, we
have ${\rm Ran\,} \Gamma_i \subset {\cal K}$. The density of ${\rm
Ran\,} L_0^* \supset {\rm Ran\,}L$ in $\cal H$ implies ${\rm
clos\,} {\rm Ran\,} \Gamma_2 = {\cal K}$. Thus, ${\rm clos\,}
\vee_{i=1,2}{\rm Ran\,} \Gamma_i = {\cal K}$ does hold.
\medskip

\noindent {\bf 2.}\,\,\,Recall the Vishik decomposition
\cite{Vishik}
\begin{equation}\label{Vishik general}
{\rm Dom\,}L_0^*={\rm Dom\,}L_0 \overset{.}+L^{-1}{\cal
K}\overset{.}+{\cal K}
\end{equation} (direct sums). By this, for a $u, v \in {\rm Dom\,}L_0^*$
one can represent
\begin{equation}\label{Vishik concrete}
u=u_0+L^{-1}g_u+h_u\,,\quad v=v_0+L^{-1}g_v+h_v \end{equation}
with $u_0, v_0 \in {\rm Dom\,}L_0$ and $g_u, h_u, g_v, h_v, \in
{\cal K}$. Therefore, with regard to $L_0\subset L_0^*$ and
$L^{-1}L_0= \mathbb I$ one has
\begin{align*}
& (L_0^*u,v)-(u,L_0^*v)=\\
& (L_0 u_0+g_u,
v_0+L^{-1}g_v+h_v)-(u_0+L^{-1}g_u+h_u, L_0v+g_v)=\\
& (L_0 u_0,v_0)+(L_0 u_0, L^{-1}g_v)+(L_0
u_0,h_v)+(g_u,v_0)+(g_u,L^{-1}g_v)+(g_u,h_v)\\
-&(u_0,L_0v_0)-(u_0,
g_v)-(L^{-1}g_u,L_0v_0)-(L^{-1}g_u,g_v)-(h_u,L_0v_0)-(h_u,g_v)\,.
\end{align*}
Numbering  by [...] the terms in the r.h.s. of the latter equality
from the beginning, we have
\begin{align*}
& [1]+[7]=0,\quad [3]=(u_0,L_0^*h_v)=0, \quad
[11]=(L_0^*h_u,v_0)=0,\\
& [2]=(L^{-1}L_0u_0,g_v)=(u_0,g_v), \quad [9]=(g_u,L^{-1}L_0
v_0)=(g_u,v_0),
\end{align*}
and arrive at
\begin{align}\label{**}
\notag & (L_0^*u,v)-(u,L_0^*v)=\\
\notag & (u_0,g_v)+(g_u,v_0)+ (g_u, L^{-1}g_v)+(g_u,h_v)\\
\notag -&(u_0,g_v)-(g_u,v_0)-(L^{-1}g_u,g_v)-(h_u,g_v)=\\
& (g_u,h_v)-(h_u,g_v)\,.
\end{align}
\medskip

\noindent {\bf 3.}\,\,\,Take a $y \in {\rm Dom\,}L_0^*$ and
represent by (\ref{Vishik concrete})
\begin{equation}\label{5*}
y\,=\,y_0+L^{-1}g_y+h_y\,.
\end{equation}
Let us check that
\begin{equation}\label{6*}
y_0=L^{-1}P_\bot L_0^* y, \quad g_y=PL_0^*y, \quad
h_y=y-L^{-1}L_0^*y\,,
\end{equation}
where $P_\bot={\mathbb I}-P$ is the projection onto $\cal H
\ominus {\cal K}$.

Indeed, since $L_0^{-1}$ is bounded on ${\rm Ran\,}L_0$, we have
${\rm Ran\,} L_0= {\rm clos\,} {\rm Ran\,}L_0$. Hence, ${\cal
H}={\rm Ran}L_0 \oplus {\cal K}$. Therefore, $P_\bot L_0^* y \in
{\rm Ran\,}L_0$ and $y=L^{-1}P_\bot L_0^*y=L_0^{-1}P_\bot L_0^*y$,
so that $y_0 \in {\rm Dom\,}L_0$. The inclusion $g_y=PL_0^*y \in
{\cal K}$ is evident. The relations $L_0^* h_y= L_0^*y-
L_0^*L^{-1}L_0^*y=L_0^*y-{\mathbb I}L_0^*y=0$ show that $h_y \in
{\cal K}$. Thus, the summands in the r.h.s. of (\ref{5*}) belong
to ${\rm Dom\,}L_0,\,L^{-1}{\cal K}$, and ${\cal K}$ respectively.

In the mean time, we have
\begin{align*}
& y=y_0+L^{-1}g_y+h_y=\langle\,{\rm see\,} (\ref{6*})\rangle\,=
L^{-1}P_\bot L_0^* y + L^{-1}PL_0^*y + y-L^{-1}L_0^*y=\\
& L^{-1}L_0^*y + y-L^{-1}L_0^*y=y\,,
\end{align*}
so that (\ref{5*}) is valid.
\medskip

\noindent {\bf 4.}\,\,\,Return to (\ref{**}). By (\ref{5*}) and
(\ref{6*}), we have
\begin{align*}
& (L_0^*u,v)-(u,L_0^*v)=\left(\left[L^{-1}L_0^*-{\mathbb I}\right]
u, P L_0^*v\right)-\left(P L_0^*u, \left[L^{-1}L_0^*-{\mathbb
I}\right]v\right)=\\
& (\Gamma_1u, \Gamma_2 v)_{\cal K} - (\Gamma_2u, \Gamma_1 v)_{\cal
K}
\end{align*}
that proves the lemma.\,\,\,\,\,\,$\blacksquare$
\medskip

Thus, the operator $L_0$ determines the Green system ${{{\mathfrak
G}{\mathfrak r}}_{L_0}}$ in a canonical way.

\subsection{System $\alpha_{L_0}$}

In its turn, ${{{\mathfrak G}{\mathfrak r}}_{L_0}}$ determines an
evolutionary dynamical system of the form
\begin{align}
\label{alpha1} & u_{tt}+L_0^*u = 0  && {\rm in}\,\,\,{\cal H}, \,\,\,t>0\\
\label{alpha2} & u|_{t=0}=u_t|_{t=0}=0 && {\rm in}\,\,\,{\cal H}\\
\label{alpha3}& \Gamma_1 u = h && {\rm in}\,\,\,{\cal K},\,\,\,t
\geqslant 0,
\end{align}
where $h=h(t)$ is a ${\cal K}$-valued function of time, $u=u^h(t)$
is a solution.

In control theory, problem (\ref{alpha1})--(\ref{alpha3}) is
referred to as a {\it dynamical system with boundary control}
(DSBC), $h$ is a {\it boundary control}, $u^h(\cdot)$ is a {\it
trajectory}, $u^h(t)$ is a {\it state} at the moment $t$. As is
clear, the system (\ref{alpha1})--(\ref{alpha3}) is determined by
the operator $L_0$ and we denote it by ${\alpha_{L_0}}$.

Recall that $L$ is the Friedrichs extension of $L_0$. Let
$L^{\frac{1}{2}}$ be the positive square root of $L$. Denote
$(\cdot)':=\frac{d}{dt}$.

For a control $h \in C^2_{\rm loc}\left([0,\infty); {\cal
K}\right)$ provided $h(0)=h'(0)=0$, an $\cal H$-valued function
\begin{equation}\label{weak solution u^f}
u^h(t)\,:=\,-h(t)+\int_0^t
L^{-\frac{1}{2}}\,\sin\left[(t-s)L^{\frac{1}{2}}\right]\,h''(s)\,ds\,,
\qquad t \geqslant 0
\end{equation}
is said to be a {\it weak solution} to
(\ref{alpha1})--(\ref{alpha3}). This definition is motivated by
the following fact. Introduce a class of smooth controls
$${\cal M}\,:=\,\{h \in C^3_{\rm loc}\left([0,\infty); {\cal K}\right)\,|\,\,h(0)=h'(0)=h''(0)=0\}.$$

\begin{Lemma}\label{Lemma 2}
If $h \in \cal M$ then $u^h$ is a classical solution to
$(\ref{alpha1})\!-\!(\ref{alpha3})$.
\end{Lemma}
{\bf Proof}
\smallskip

\noindent {\bf 1.}\,\,Assuming $h \in \cal M$, let us derive a
relevant representation for the weak solution.

Take a $y \in \cal H$. Representing $L=\int_0^\infty
\lambda\,dE_\lambda$ via the spectral measure $E_\lambda$ of $L$
and integrating by parts, one has
\begin{align*}
& \frac{d}{ds}\,\int_0^\infty \frac{\cos
\sqrt{\lambda}(t-s)}{\lambda}\,d\left(E_\lambda
h''(s),y\right)=\frac{d}{ds}\,\left(h''(s),
L^{-1}\,\cos\left[(t-s)L^{\frac{1}{2}}\right]y\right)=\\
& \left(h'''(s),
L^{-1}\,\cos\left[(t-s)L^{\frac{1}{2}}\right]y\right)+\left(h''(s),
L^{-\frac{1}{2}}\,\sin\left[(t-s)L^{\frac{1}{2}}\right]y\right)=\\
& \left(L^{-1}\,\cos\left[(t-s)L^{\frac{1}{2}}\right]h'''(s),
y\right)+\left(L^{-\frac{1}{2}}\,\sin\left[(t-s)L^{\frac{1}{2}}\right]h''(s),
y\right).
\end{align*}
Applying $\int_0^t(\dots)\,ds$, we get
\begin{align*}
& \int_0^\infty\frac{1}{\lambda}\,d \left(E_\lambda
h'''(t),y\right)= \left(L^{-1}h'''(t), y\right)=\\
& \left(\int_0^t
L^{-1}\,\cos\left[(t-s)L^{\frac{1}{2}}\right]h'''(s),
y\right)+\left(\int_0^t
L^{-\frac{1}{2}}\,\sin\left[(t-s)L^{\frac{1}{2}}\right]h''(s),
y\right).
\end{align*}
By arbitrariness of $y$, we get
\begin{align*}
\int_0^t
L^{-\frac{1}{2}}\,\sin\left[(t-s)L^{\frac{1}{2}}\right]h''(s)\,ds=
L^{-1}\int_0^t\left\{{\mathbb I}-
\cos\left[(t-s)L^{\frac{1}{2}}\right]\right\}h'''(s)\,ds.
\end{align*}
Therefore, for $h \in \cal M$ one can write (\ref{weak solution
u^f}) in the form
\begin{align}\label{12*}
u^h(t)=-h(t)+ L^{-1}\int_0^t\left\{{\mathbb I}-
\cos\left[(t-s)L^{\frac{1}{2}}\right]\right\}h'''(s)\,ds
\end{align}
with $h \in {\cal K}$ and $L^{-1}\int_0^t(\dots)\in {\rm Dom\,}
\subset {\rm Dom\,}L_0^*$. Hence, we have $u^h(t) \in {\rm
Dom\,}L_0^*$ for all $t \geqslant 0$.
\medskip

\noindent {\bf 2.}\,\,\,Show that (\ref{12*}) provides a classical
solution to (\ref{alpha1})--(\ref{alpha3}).

Differentiation in (\ref{12*}) implies
\begin{align}\label{13*}
& u^h_t(t)=-h'(t)+ L^{-\frac{1}{2}}\int_0^t
\sin\left[(t-s)L^{\frac{1}{2}}\right]h'''(s)\,ds, \\
& u^h_{tt}(t)=-h''(t)+ \int_0^t
\cos\left[(t-s)L^{\frac{1}{2}}\right]h'''(s)\,ds.
\end{align}
Therefore,
\begin{align*}
& u^h_{tt}(t)+L_0^* u^h(t)=-h''(t)+ \int_0^t
\cos\left[(t-s)L^{\frac{1}{2}}\right]h'''(s)\,ds+\\
& \int_0^t\left\{{\mathbb I}-
\cos\left[(t-s)L^{\frac{1}{2}}\right]\right\}h'''(s)\,ds=
-h''(t)+\int_0^th'''(s)\,ds\,=\,0,
\end{align*}
so that (\ref{alpha1}) does hold.
\smallskip

As is seen from (\ref{12*}),(\ref{13*}), $u^h(0)=u^h_t(0)=0$,
i.e., the initial conditions (\ref{alpha2}) are fulfilled.
\smallskip

Applying $\Gamma_1$ in (\ref{12*}), we have
\begin{align*}
& \Gamma_1 u^h(t)=\left(L^{-1}L_0^*-{\mathbb I}\right)\left(-h(t)+
L^{-1}\int_0^t\left\{{\mathbb I}-
\cos\left[(t-s)L^{\frac{1}{2}}\right]\right\}h'''(s)\,ds\right)=\\
& L^{-1}\int_0^t\left\{ {\mathbb I}-
\cos\left[(t-s)L^{\frac{1}{2}}\right]\right\}h'''(s)\,ds+h(t)\,-\\
& L^{-1}\int_0^t\left\{{\mathbb
I}-\cos\left[(t-s)L^{\frac{1}{2}}\right]\right\}h'''(s)\,ds\,=\,h(t).
\end{align*}
Hence, the `boundary condition' (\ref{alpha3}) is
fulfilled.\,\,\,$\blacksquare$
\smallskip

Note in addition, that one can prove a uniqueness of  the
classical $u^h$.

\subsection{Controllability}
A set of all possible states of the system ${\alpha_{L_0}}$
$${\cal U}^t:=\{u^h(t)\,|\,\, u^h \,\,{\rm is\,\,a\, week\,\,solution\,\,to \,\,
(\ref{alpha1})\!-\!(\ref{alpha3})}\}$$ is said to be {\it
reachable} (at the moment $t \geqslant 0$). Representing
(\ref{weak solution u^f}) in the convolution form
\begin{equation}\label{u^h is a convolution}
u^h(t)\,=\,\int_0^t\left\{-(t-s){\mathbb I}+
L^{-\frac{1}{2}}\,\sin\left[(t-s)L^{\frac{1}{2}}\right]\right
\}\,h''(s)\,ds\,,
\end{equation}
one can easily see that ${\cal U}^t$ extends as $t$ grows.

Also, define a {\it total} reachable set
$${\cal U}\,:=\,\bigvee_{t \geqslant 0}{\cal U}^t$$ and a {\it defect
subspace}$${\cal D}\,:=\,{\cal H} \ominus {\rm clos\,} \cal U\,.$$
The system ${\alpha_{L_0}}$ is said to be {\it controllable}, if
the relation
\begin{equation}\label{contr}
{\rm clos\,} {\cal U}\,=\,\cal H
\end{equation}
is valid or, equivalently, if ${\cal D}=\{0\}$.

Below we establish the necessary and sufficient conditions on the
operator $L_0$, which provide controllability of the system
${\alpha_{L_0}}$. Taking into account the well-known similar
results and general principles of system theory \cite{KFA}, these
conditions are quite expectable: $L_0$ has to be a {\it completely
non-self-adjoint} (c.n.s.a.) operator.

Recall the definitions. We say that a symmetric operator $A$ has a
self-adjoint part in a (nonzero) subspace ${\cal N} \subset \cal
H$ if
\begin{itemize}
\item the lineal set ${\cal N} \cap {\rm Dom\,}A$ is dense in
$\cal N$ \item the embedding $A \left[{\cal N} \cap {\rm
Dom\,}A\right]\,\subset \cal N$ holds \item the operator
$A|_{{\cal N} \cap {\rm Dom\,}A}$ is self-adjoint in $\cal N$.
\end{itemize}
A symmetric operator $A$ is said to be c.n.s.a. if it has a
self-adjoint part in no nonzero subspace in $\cal H$.
\begin{Theorem}\label{Theorem 1}
The system ${\alpha_{L_0}}$ is controllable if and only if $L_0$
is a c.n.s.a. operator.
\end{Theorem}
{\bf Proof}
\medskip

\noindent{\bf Necessity}\,\,\, Let $L_0$ have a self-adjoint part
in ${\cal N} \subset \cal H$. Fix an $h \in {\cal K}$. Take a $g
\in \cal N$ and represent $g=L_0 \widetilde g$ with $\widetilde g
\in {\cal N} \cap {\rm Dom\,}L_0$. The latter is possible because
$L_0|_{\cal N}$ is a positive definite boundedly invertible
operator in $\cal N$. In view of
$$(g,h)=(L_0 \widetilde g,h)=(\widetilde g, L_0^*)=0$$
we have ${\cal N} \perp {\cal K}$, i.e.,
\begin{equation}\label{1111}
{\cal K} \subset {\cal N}^\bot={\cal H} \ominus {\cal N}
\end{equation} holds.
\smallskip

Recall that $L$ is the Friedrichs extension: $L_0 \subset L
\subset L_0^*$. For a $g \in \cal N$ we have $$L^{-1}g=L^{-1}L_0
\widetilde g = \widetilde g \in \cal N$$ that implies $L^{-1}{\cal
N} \subset \cal N$. Since $L^{-1}$ is self-adjoint, we have
$L^{-1}{\cal N}^\bot \subset {\cal N}^\bot$, i.e., $L^{-1}$ is
reduced by $\cal N$. The latter leads to
\begin{equation}\label{14*}
L^{-\frac{1}{2}}{\cal N} \subset {\cal N}, \qquad
L^{-\frac{1}{2}}{\cal N}^\bot \subset {\cal N}^\bot.
\end{equation}

By (\ref{1111}), in the r.h.s. of (\ref{weak solution u^f}), one
has $h(t),h''(s) \in {\cal K} \subset {\cal N}^\bot$. By
(\ref{14*}), the integral in (\ref{weak solution u^f}) belongs to
${\cal N}^\bot$. As a result, $u^h(t) \in {\cal N}^\bot$ holds for
all $t \geqslant 0$, i.e., trajectories $u^h$ of the system
${\alpha_{L_0}}$ do not leave the subspace ${\cal N}^\bot$.
Therefore, ${\cal U} \subset {\cal N}^\bot$ that leads to ${\cal
D} \not=\{0\}$. So, the system ${\alpha_{L_0}}$ is not
controllable.
\smallskip

Thus, if ${\alpha_{L_0}}$ is controllable then $L_0$ is c.n.s.a.
\bigskip

{\bf Sufficiency}\,\,\,Assume that ${\alpha_{L_0}}$ is not
controllable, i.e., ${\cal D}\not=\{0\}$. It will be shown that
$L_0$ has to have a self-adjoint part in $\cal D$ and, hence, is
not a c.n.s.a. operator.
\smallskip

\noindent{\bf 1.}\,\,\,Take a nonzero $y \in \cal D$. For any
(admissible) $h \in C^2_{\rm loc}\left([0,\infty); {\cal
K}\right)$ and $t>0$, we have \begin{align}\label{15*} \notag &
0=(y,u^h(t))=\langle {\rm see\,}(\ref{weak solution u^f}),
(\ref{u^h is a convolution})\rangle=\\
\notag & \left(y, \int_0^t\left\{-(t-s)h''(s)+
L^{-\frac{1}{2}}\sin\left[(t-s)L^{\frac{1}{2}}\right]\right
\}\,h''(s)\,ds\right)=\\
& \int_0^t \left(-(t-s)y + w^y(t-s), h''(s)\right)\,ds,
\end{align}
where
$$w^y(\eta)\,:=\,L^{-\frac{1}{2}}\,\sin\left[\eta L^{\frac{1}{2}}\right] y,
\qquad \eta \geqslant 0.$$

Fix a $k \in {\cal K}$. Choose a sequence of controls
$h_j(s)=\varphi_j(s) k$, where $\varphi_j \in C^\infty_0(0,t)$ are
such that $\varphi''_j(s) \to \delta(s)$ (the Dirac
delta-function) as $j \to \infty$. For such $h_j(\cdot)$, the
limit passage in (\ref{15*}) implies $$0=(-ty + w^y(t), k).$$ By
arbitrariness of $k$, we conclude that
\begin{equation}\label{16*}
-ty + w^y(t) \in {\cal K}^\bot, \qquad t \geqslant 0.
\end{equation}
Converting these considerations, we easily obtain the following
result.
\begin{Proposition}\label{Prpopsition 1}
The embedding $y \in \cal D$ holds if and only if $(\ref{16*})$ is
valid.
\end{Proposition}
Denote $-ty + w^y(t)=:p(t) \in {\cal K}^\bot$ and represent
$$y=\frac{1}{t}\,L^{-\frac{1}{2}}\,\sin\left[tL^{\frac{1}{2}}\right]
-\frac{p(t)}{t}\,.$$ The operator
$L^{-\frac{1}{2}}\,\sin\left[tL^{\frac{1}{2}}\right]$ is bounded.
By the latter, tending $t \to \infty$ we get $y=-\lim t^{-1}p(t)
\in {\cal K}^\bot$. Returning to (\ref{16*}), we conclude that
\begin{equation}\label{17*}
y,\,w^y(t) \in {\cal K}^\bot, \qquad t \geqslant 0.
\end{equation}
\medskip

\noindent {\bf 2.}\,\,\,Consider an auxiliary dynamical system
\begin{align}
\label{18*}& w_{tt}+Lw=0 && {\rm in}\,\,\,{\cal H}, \quad t>0\\
\label{19*}& w|_{t=0}=0, \,\,\,w_t|_{t=0}=y&& {\rm in}\,\,\,\cal H
\end{align}
with $y \in \cal D$ chosen above. Its solution is of the
well-known form
$$w^y(t)=L^{-\frac{1}{2}}\,\sin\left[t L^{\frac{1}{2}}\right]
y$$ (see, e.g., \cite{BirSol}).

Denote $J:=\int_0^t(\cdot)\,ds$. Applying $J^2$ in (\ref{18*})
with regard to (\ref{19*}), we have $$w^y(t)-ty\,=\,- L\left(J^2
w^y\right)(t).$$ This implies
\begin{equation}\label{20*}
L^{-1}\left(-ty+w^y(t)\right)\,=\,- \left(J^2 w^y\right)(t) \in
{\cal K}^\bot
\end{equation}
in view of (\ref{17*}). In the mean time, we have
$$L^{-1}w^y(t)= L^{-1}L^{-\frac{1}{2}}\,\sin\left[t L^{\frac{1}{2}}\right]
y = L^{-\frac{1}{2}}\,\sin\left[t
L^{\frac{1}{2}}\right]\left(L^{-1} y\right)= w^{L^{-1}y}(t)\,.$$
Hence, (\ref{20*}) implies
\begin{equation*}
-tL^{-1}y+w^{L^{-1}y}(t) \in {\cal K}^\bot, \qquad t \geqslant 0.
\end{equation*}
In accordance with Proposition \ref{Prpopsition 1}, the latter is
equivalent to $L^{-1}y \in \cal D$.

Thus, beginning with $y \in \cal D$, we arrive at $L^{-1}y \in
\cal D$, i.e., the defect space reduces the operator $L^{-1}$:
$$L^{-1}\cal D\,\subset \cal D\,.$$
\medskip

\noindent {\bf 3.}\,\,\,The part $L^{-1}|_{\cal D}$ is a
self-adjoint invertible operator in $\cal D$. Therefore,
$L^{-1}\cal D$ is dense in $\cal D$. Hence, the operator $L$ has
the part $L|_{\cal D}$, which is a (densely defined) self-adjoint
operator in $\cal D$.

Show that $L|_{\cal D}=L_0$. Indeed, by (\ref{17*}) one has $$\cal
D \subset {\cal K}^\bot=\cal H \ominus {\cal K}={\rm clos\,} {\rm
Ran\,} L_0 = {\rm Ran\,} L_0\,.$$ Therefore, by $L_0 \subset L$ we
have
$${\rm Dom\,}L|_{\cal D}=L^{-1}{\cal D} \subset L^{-1}{\rm Ran\,}L_0 =
L_0^{-1}{\rm Ran\,} L_0 = {\rm Dom\,}L_0,$$ and $L|_{\cal D} x =
L_0 x$ for all $x \in L^{-1}\cal D$.
\medskip

Thus, if the system ${\alpha_{L_0}}$ is not controllable then
$L_0$ has a self-adjoint part in $\cal D$, i.e., is not a c.n.s.a.
operator. \,\,\,$\blacksquare$

\section{Wave spectrum}

\subsection{Inflation}
We use the term `lattice' in its general meaning \cite{Birkhoff}:
a {\it lattice} is a partially ordered set provided with the
operations $x \vee y={\rm sup\,}\{x,y\},\,\,x \wedge y ={\rm
inf\,}\{x,y\}$. However, we deal with the concrete lattices
endowed with additional structures (complement, topology, ets).
\smallskip

\noindent {\bf Definition}\,\,\,Let ${{\mathfrak L}(\cal H)}$ be
the lattice of the (closed) subspaces of $\cal H$ with the partial
order $\subseteq$ and operations ${\cal A} \wedge {\cal B}={\cal
A} \cap {\cal B},\,\,{\cal A} \vee {\cal B}={\rm
clos\,}\{a+b\,|\,\,a \in {\cal A},\, b \in \cal B\},\,\,{\cal A}
\mapsto {\cal A}^\bot$. Also, it possesses the least and greatest
elements $\{0\}$ and $\cal H$.

By $P_{\cal A}$ we denote the (orthogonal) projection in $\cal H$
onto $\cal A$. Topology on ${{\mathfrak L}(\cal H)}$ is determined
by convergence of the projections on the corresponding subspaces.
Namely, we write $\cal A_j \to \cal A$ if $s\!-\!\lim P_{{\cal
A}_j}=P_{\cal A}$.
\medskip

An {\it inflation} is a family of maps $I=\{I^t\}_{t \geqslant
0},\,\,I^t: {{\mathfrak L}(\cal H)} \to {{\mathfrak L}(\cal H)}$
with the properties
\begin{itemize}
\item $I^0\,=\,{\rm id},\,\,\,I^t\{0\}\,=\,\{0\}$ \item ${\cal A}
\subseteq \cal B$ and $s \leqslant t$ imply $I^s {\cal A}
\subseteq I^t \cal B$.
\end{itemize}
\medskip

\noindent{\bf Inflation ${I_{L_0}}$}\,\,\, As is shown in
\cite{BArXiv}, with each operator $L_0$ of the class under
consideration one associates an inflation $I_{L_0}$ in the
following way. Fix a subspace $\cal A \in {{\mathfrak L}(\cal H)}$
and consider a dynamical system
\begin{align*}
& v_{tt}+L v = a  && {\rm in}\,\,\,{\cal H}, \,\,\,t>0\\
& v|_{t=0}=v_t|_{t=0}=0 && {\rm in}\,\,\,{\cal H}\\
\end{align*}
where $a=a(t)$ is an $\cal A$-valued function of time, $v=v^a(t)$
is a solution. By the well-known Duhamel formula, one has
$$v^a(t)=\int_0^t
L^{-\frac{1}{2}}\,\sin\left[(t-s)L^{\frac{1}{2}}\right]a(s)\,ds,
\qquad t \geqslant 0.$$ The reachable sets of the system are
$${\cal V}^t_{\cal A}:=\left\{v^a(t)\,|\,\,a \in L_2^{\rm loc}\left([0,\infty); {\cal A}\right)\right\}.$$
As is easy to see, ${\cal V}^t_{\cal A}$ extends as ${\cal A}$
extends and/or $t$ grows. Define the family
${I_{L_0}}:=\{I_{L_0}^t\}_{t \geqslant 0}$ by
$$I_{L_0}^t{\cal A}:={\rm clos\,} {\cal V}^t_{\cal A}\,\,\,\,\,{\rm for}\,\,t>0, \qquad I_{L_0}^0:={\rm id}.$$
\begin{Proposition}
The family ${I_{L_0}}$ is an inflation.
\end{Proposition}
Proof see in \cite{BArXiv}.

Actually, ${I_{L_0}}$ is determined not by $L_0$ but its extension
$L$. As such, an inflation is well-defined for any bounded from
bellow self-adjoint operator.

\subsection{Set ${\Omega_{L_0}}$}
{\bf Lattice ${{\mathfrak L}_{L_0}}$}\,\,\,Recall that a {\it
sublattice} is a subset of ${{\mathfrak L}(\cal H)}$, which is
invariant with respect to all the lattice operations. Each
sublattice necessarily contains $\{0\}$ and $\cal H$.

We say that a (sub)lattice ${\mathfrak L} \subset {{\mathfrak
L}(\cal H)}$ is invariant w.r.t. an inflation $I$ if $I^t
{\mathfrak L} \subseteq {\mathfrak L}$ holds for $t \geqslant 0$.
\medskip

By ${{\mathfrak L}_{L_0}}$ we denote the {\it minimal sublattice}
in ${\mathfrak L}$, which
\begin{itemize} \item
contains all reachable subspaces ${\rm clos\,}{\cal U}^t,\,t
\geqslant 0$ \item is invariant with respect to the inflation
${I_{L_0}}$ \item is closed in the above-mentioned topology on
${{\mathfrak L}(\cal H)}$.
\end{itemize}

\noindent{\bf Lattice $\overline {{I_{L_0}}{{\mathfrak
L}_{L_0}}}$}\,\,\, Let $\cal F$ be a set of ${{\mathfrak L}(\cal
H)}$-valued functions of $t \geqslant 0$. This set is also a
lattice w.r.t. the point-wise order, operations, and topology:
\begin{align*}
& \{f \leqslant g \} \Longleftrightarrow \{f(t) \subseteq
g(t),\,\,t \geqslant 0\}, \quad (f \vee g)(t):=f(t)\vee g(t),\\
& (f \wedge g)(t):=f(t)\wedge g(t),
\,\,(f^\bot)(t):=(f(t))^\bot,\,\,(\lim f_j)(t):=\lim (f_j(t))\,.
\end{align*}
The least and greatest elements of $\cal F$ are the functions
equal $\{0\}$ and $\cal H$ identically. We denote them by $0_{\cal
F}$ and $1_{\cal F}$ respectively.

An inflation $I$ can be regarded as a map from ${{\mathfrak
L}(\cal H)}$ to $\cal F$ acting by the rule $(I{\cal
A})(t):=I^t{\cal A},\,\,t \geqslant 0$. If ${\mathfrak L}$ is
invariant w.r.t. $I$ then the image $I{\mathfrak L}$, as well as
its closure $\overline {I{\mathfrak L}}$ are sublattices in $\cal
F$. Both of them contain $0_{\cal F}$.
\smallskip

The operator $L_0$ determines the lattice $\overline
{{I_{L_0}}{{\mathfrak L}_{L_0}}}$.
\medskip

\noindent{\bf Atoms}\,\,\,Let $\cal P$ be a partially ordered set
with the least element $0$. An element $a \in \cal P$ is said to
be an {\it atom} if $a \not= 0$ and $b\leqslant a$ implies $b=a$
\cite{Birkhoff}. By ${\rm At} \cal P$ we denote the set of atoms.
\smallskip

The key object of the paper is the set
$${\Omega_{L_0}}\,:=\,{\rm At}\overline {{I_{L_0}}{{\mathfrak L}_{L_0}}}$$ that we call
a {\it wave spectrum} of the operator $L_0$.
\medskip

\noindent{\bf Remark}\,\,\,Certain additional assumptions on $L_0$
provide ${\Omega_{L_0}}\not=\emptyset$ \cite{BArXiv}. There is
$L_0$ such that its wave spectrum consists of a single point. A
conjecture is that ${\Omega_{L_0}}\not=\emptyset$ does hold ever.

\subsection{Space $\left({\Omega_{L_0}},\beta\right)$}
Here the wave spectrum is endowed with relevant structures.
\medskip

\noindent{\bf Topology}\,\,\,Recall that atoms are ${{\mathfrak
L}_{L_0}}$-valued functions of time. Fix an atom $a \in
{\Omega_{L_0}}$. A set
$$B_r[a]:=\{b
\in {\Omega_{L_0}}\,|\,\,\exists t>0\,\,\,{\rm
s.t.}\,\,\{0\}\not=b(t)\leqslant a(r)\}, \qquad r>0$$ is said to
be a {\it ball}, $a$ and $r$ are its  center and radius.
\begin{Lemma}\label{Lemma 3}
The system of balls $\{B_r[a]\,|\,\,a \in
{\Omega_{L_0}},\,\,r>0\}$ is a base of topology.
\end{Lemma}
{\bf Proof}\,\,\,One has to check the characteristic properties of
a base:
\begin{enumerate}
\item for any $a \in {\Omega_{L_0}}$, there is a ball $B \ni a$
\item for an atom $a \in {\Omega_{L_0}}$ and the balls $B_1,\,B_2$
such that $a \in B_1 \cap B_2$, there is a ball $B$ such that $a
\in B \subset B_1 \cap B_2$
\end{enumerate}
(see, e.g., \cite{Kelley}).
\medskip

\noindent{\bf 1.}\,\,\,Take an $a=a(\cdot)\in {\Omega_{L_0}}$. For
any $r > t_0:={\rm inf}\left\{t>0\,|\,\,a(t)\not= \{0\}\right\}$
and $t \in (t_0, r]$, one has $\{0\}\not=a(t) \leq a(r)$, i.e., $a
\in B_r[a]$.
\medskip

\noindent{\bf 2.}\,\,\,Let $a \in B^{r_1}[a_1] \cap B^{r_2}[a_2]$,
so that both of $B^{r_i}[a_i]$ are nonempty. Choose $t_i$ such
that $\{0\}\not=a(t_i) \leqslant a_i(r_i)$ and denote $r:={\rm
min}\{t_1, t_2\}$. By the choice, one has $\{0\}\not=a(r)
\leqslant a_i(r_i)$ that implies $B_r[a]\not=\emptyset$.

For any $b \in B_r[a]$, there is a $t>0$ such that
$\{0\}\not=b(t)\leqslant a(r) \leqslant a_i(r_i)$. By the latter
inequality, one has $b \in B_{r_i}[a_i]$. Hence $B_r[a] \subset
B^{r_1}[a_1] \cap B^{r_2}[a_2]$.\,\,\,$\blacksquare$
\medskip

The base $\{B_r[a]\,|\,\,a \in {\Omega_{L_0}},\,\,r>0\}$
determines the (unique) topology, in which an open set is a sum of
balls \cite{Kelley}. We call it a {\it ball topology} and denote
by $\beta$. So, we get a topological space $\left({\Omega_{L_0}},
\beta\right)$.
\medskip

\noindent{\bf Boundary}\,\,\,Return to the DSBC ${\alpha_{L_0}}$.
The family of reachable subspaces
$${\mathfrak u}_{L_0}\,:=\,\left\{{\rm clos\,} \cal U^t\right\}_{t \geqslant 0}$$
can be regarded as an ${{\mathfrak L}_{L_0}}$-valued function of
time. As such, ${\mathfrak u}_{L_0}$ is an element of the lattice
$\overline {{I_{L_0}}{{\mathfrak L}_{L_0}}} \subset \cal F$ and
can be compared with its atoms. Thus, the set
$$\partial {\Omega_{L_0}}\,:=\,\left\{a \in {\Omega_{L_0}}\,|\,\,\,a \leqslant {\mathfrak u}_{L_0}\right\}$$
is well defined and said to be a {\it boundary} of the wave
spectrum.

\section{Illustration}
\subsection{Manifold}
Let $\Omega$ be a $C^\infty$-smooth compact Riemannian manifold of
dimension $n \geqslant 2$ with the boundary $\partial \Omega$, $g$
the metric tensor, $-\Delta$ the scalar Beltrami-Laplace operator.
Recall that in local coordinates one has $$-\Delta=-[{\rm
det\,}g]^{-\frac{1}{2}}\frac{\partial}{\partial x_i}[{\rm
det\,}g]^{\frac{1}{2}}g^{ij}\frac{\partial}{\partial x_j}\,.$$ By
$\nu$ we denote an outward normal to $\partial \Omega$;
$\partial_\nu$ is differentiation w.r.t. the normal.

The manifold is endowed with volume form $dv$, so that the (real)
Hilbert space ${\cal H}:=L_2(\Omega)$ with the inner product
$$(u,v)=\int_\Omega u\,v\,dv$$ is well defined.

The boundary $\partial \Omega$ is endowed with canonical (induced
by the tensor $g|_{\partial \Omega}$) metric and volume (surface)
element $d\sigma$. In the space ${\cal B}:=L_2(\partial \Omega)$,
the inner product is
$$(f,g)_{\cal B}=\int_{\partial \Omega} f g\,d\sigma.$$

\subsection{Operators}
Our basic operator is the {\it minimal Laplacian} $L_0: {\cal H}
\to {\cal H},\,\,{\rm Dom\,}L_0=\{y \in
H^2(\Omega)\,|\,\,y={\partial_\nu y}=0\,\,{\rm on}\,\,\,\partial
\Omega\},\,\,L_0=-\Delta y$, where $H^2(\Omega)$ is the Sobolev
class. Operator $L_0$ is positive definite and symmetric, its
defect indexes are $n_\pm=\infty$.

The operator $L_0^*$ is the {\it maximal Laplacian}, which is
defined on ${\rm Dom\,}L_0^*=\{y \in {\cal H}\,|\,\,\Delta y \in
{\cal H}\}$ (here $\Delta y$ is understood in the sense of
distributions) and acts by $L_0^*y=-\Delta y$. Its null subspace
consists of harmonic functions: $${\cal K}\,=\,{\rm Ker\,}L_0^*
=\{y \in {\cal H}\,|\,\,\Delta y=0\,\,\,{\rm in}\,\,\,\,\Omega
\backslash
\partial \Omega\}\,.$$ Also, we use the notation
$L_0^*=-\Delta_{\rm max}$.

The {\it Friedrichs extension} $L \supset L_0$ is defined on ${\rm
Dom\,}L= \{y \in H^2(\Omega)\,|\,\,y=0\,\,{\rm on}\,\,\,\partial
\Omega\}$ and acts by $Ly=-\Delta y$.
\smallskip

One more operator associated with the manifold is the {\it
harmonic continuation} $\Pi: \cal B \to \cal H$ defined by the
relations $\Delta \Pi f = 0$ in $\Omega$, and $\left(\Pi
f\right)|_{\partial \Omega}=f$. As is well-known in elliptic PDE
theory, $\Pi$ is a compact injective operator; its adjoint
$\Pi^*:\cal H \to \cal B$ is
$$\Pi^*\,=\,{\partial_\nu}\,L^{-1}\,.$$

\subsection{Green system}
As one can check, for a $y \in {\rm Dom\,}L_0^*$, the
decomposition (\ref{Vishik concrete}) is $$y=y_0+L^{-1}g_y+h_y$$
with
\begin{equation}\label{21*}
h_y\,=\,\Pi\left(y|_{\partial \Omega}\right), \quad g_y\,=\,
\left(\Pi^*\right)^{-1}\left[{\partial_\nu
y}-\Lambda\left(y|_{\partial \Omega}\right)\right],
\end{equation}
where $\Lambda: \cal B \to \cal B$ is the {\it
Dirichlet-to-Neumann map} defined by $$\Lambda
f:={\partial_\nu}\,\Pi f \qquad {\rm on}\,\,\,\,\partial \Omega.$$
Note that the right hand sides in (\ref{21*}) have to be
understood properly; they are well defined for smooth enough $y$'s
and then extended on all $y \in {\rm Dom\,}L_0^*$ by relevant
continuity \cite{Vishik}, \cite{Ryzh}.

By (\ref{21*}), we have
\begin{align}\label{22*}
\notag & (L_0^*u, v)_{\cal H}-(u, L_0^*v)_{\cal H}=\langle {\rm
see\,}(\ref{**})\rangle=(g_u, h_v)_{\cal H} - (h_u, g_v)_{\cal H}=\\
\notag & \left(\left(\Pi^*\right)^{-1}\left[{\partial_\nu
u}-\Lambda\left(u|_{\partial \Omega}\right)\right]\,,\,
\Pi\left(v|_{\partial
\Omega}\right)\right)_{\cal H}-\\
\notag &\left(\Pi\left(u|_{\partial \Omega}\right)\,,
\,\left(\Pi^*\right)^{-1}\left[{\partial_\nu
v}-\Lambda\left(v|_{\partial \Omega}\right)\right]\right)_{\cal H}=\\
& (\gamma_1 u, \gamma_2 v)_{\cal B} - (\gamma_2 u, \gamma_1
v)_{\cal B}\,,
\end{align}
where $\gamma_i: \cal H \to \cal B$,
$$\gamma_1:=(\,\cdot\,)|_{\partial \Omega}, \quad \gamma_2:=
\left[{\partial_\nu}-\Lambda\right]\left[(\,\cdot\,)|_{\partial
\Omega}\right]$$ are the canonical (by Vishik) {\it boundary
operators}: see \cite{Vishik}, sec 6.

As a result, the canonical Green system associated with the
manifold is \begin{equation}\label{Green Omega}{{\mathfrak
G}{\mathfrak r}}_{\Omega}=\{L_2(\Omega), L_2(\partial \Omega);
-\Delta_{\rm max}, \gamma_1, \gamma_2\}.\end{equation}

\subsection{System $\alpha_\Omega$}
In accordance with (\ref{22*}), (\ref{Green Omega}), the system
${\alpha_{L_0}}=:\alpha_\Omega$ on the manifold takes the form
\begin{align*}
& u_{tt}-\Delta u = 0  && {\rm in}\,\,
           \left(\Omega \backslash \partial \Omega\right) \times (0,\infty)\\
& u|_{t=0}=u_t|_{t=0}=0 && {\rm in}\,\,\Omega\\
& u = f && {\rm in}\,\,\partial \Omega \times (0,\infty),
\end{align*}
where $f \in L_2^{\rm loc}\left(\partial \Omega \times
(0,\infty)\right)$ is a {\it boundary control}, $u=u^f(x,t)$ is a
solution. The solution describes a {\it wave}, which is initiated
by the boundary source and propagates into the manifold. The speed
of propagation is finite (equal $1$).

The system $\alpha_\Omega$ is controllable \cite{BArXiv}.
Moreover, for a compact $\Omega$ one has
\begin{equation*}
{\rm clos\,} {\cal U}^t\,=\,{\cal H}\,, \qquad t>\underset{x \in
\Omega}{\rm min}\,{\rm dist\,}(x,\partial \Omega)
\end{equation*}
(see \cite{BIP07}).

\subsection{Wave spectrum}
In \cite{BArXiv} a class of the so-called {\it simple manifolds}
is introduced. Roughly speaking, simplicity means that the group
of symmetries (isometries) of $\Omega$ is trivial. This property
is generic: any smooth compact manifold can be made simple by
arbitrarily small variations of its boundary.
\smallskip

As is shown in \cite{BArXiv}, if $\Omega$ is simple then there is
a canonical bijection $$\Omega_{-\Delta_{\rm min}} \ni a_x
\leftrightarrow x_a \in \Omega$$ between atoms and points, which
relates the balls and boundaries:
$$B_r[a_x]\leftrightarrow \{x' \in \Omega\,|\,\,{\rm dist\,}(x', x_a)<r\}, \quad
\partial \Omega_{-\Delta_{\rm min}} \leftrightarrow \partial \Omega\,.$$
Thus, a simple manifold is identical (isometric) to its wave
spectrum. It is the fact, which is used in inverse problems for
{\it reconstruction}.

Namely, each kind of traditional inverse data (response operator
\cite{BIP07}, Weyl function \cite{Ryzh}, characteristic function
\cite{Shtraus}) determines the operator $L_0=-\Delta_{\rm min}$ up
to unitary equivalence. By this, given the inverse data of a {\it
simple} manifold, one can determine a unitary copy $\widetilde
L_0$, find its wave spectrum $\Omega_{\widetilde L_0}$ and thus
recover the manifold up to isometry.
\medskip

Note that a {\it reconstruction up to isometry} is the most that
we can hope for. Assume that we are given with the boundary
inverse data of a certain manifold $\Omega$. Assume that another
$\Omega^\prime$ is {\it isometric} to $\Omega$ and has the same
boundary: $\partial \Omega^\prime =
\partial \Omega$. As is easy to recognize, the boundary data of
$\Omega^\prime$ and $\Omega$ are identical. Therefore, in
principle, these data do not determine the original $\Omega$
uniquely. In such a situation, the only relevant understanding of
`to recover' is to provide a {\it representative} of the class of
manifolds with the given data. The wave spectrum
$\Omega_{\widetilde L_0}$ does provide such a representative.

\end{document}